\documentclass[12pt,letterpaper]{article}
\usepackage[margin=1in]{geometry}
\usepackage[utf8]{inputenc}

\usepackage{booktabs} 

\usepackage{enumitem}

\usepackage{amsmath}
\usepackage{mathtools}
\usepackage{amsfonts}
\usepackage{amsthm}
\usepackage{amssymb}
\usepackage{bm}
\usepackage{graphics}
\usepackage{graphicx}
\usepackage{caption}
\usepackage{enumitem}
\usepackage{eurosym}
\usepackage[colorinlistoftodos]{todonotes}
\usepackage{algorithm2e}
\usepackage{subfig}
\usepackage{wrapfig}
\usepackage{url}
\usepackage{wasysym}
\usepackage{resizegather}

\newcommand{\bma}{\begin{bmatrix}}
	\newcommand{\ebma}{\end{bmatrix}}

\newcommand{\bd}{\mathbf}

\newcommand{\la}{\lambda}

\newtheorem{lemma}{Lemma}

\newtheorem{theorem}{Theorem}
\newtheorem{prop}{Proposition}

\DeclarePairedDelimiter\floor{\lfloor}{\rfloor}

\title{
Optimizing Curbside Parking Resources Subject to Congestion Constraints}

\author{Chase Dowling, Tanner Fiez, Lillian Ratliff, and Baosen Zhang
\thanks{The authors have been supported in part by NSF grants CNS-1646912 and CNS-1634136. C. Dowling was also supported in part by the Washington Clean Energy Institute.}
\thanks{C. Dowling, T. Fiez, L. Ratliff and B. Zhang are with the Department of Electrical Engineering, University of Washington,
	    Seattle, WA 98195, USA
	    Emails: \{cdowling,fiezt,ratliffl,zhangbao\}@uw.edu}}

\begin{document}

\maketitle
\thispagestyle{empty}
\pagestyle{empty}

\begin{abstract}
    To gain theoretical insight into the relationship between parking scarcity and
congestion, we describe block-faces of curbside parking as a network of queues.
Due to the nature of this network, canonical queueing network
results are not available to us. We present a new kind of queueing network
subject to customer rejection due to the lack of available servers. We provide conditions for such networks to be stable, a computationally tractable ``single node'' view of such a network, and show that maximizing the occupancy through price control of such queues, and subject to constraints on the allowable congestion between queues searching for an available server, is a convex optimization problem. We demonstrate an application of this method in the Mission District of San Francisco; our results suggest congestion due to drivers searching for parking stems from an inefficient spatial utilization of parking resources.

\end{abstract}

\section{INTRODUCTION}
Drivers in densely populated urban districts often find that desirable parking
close to their destination is unavailable or prohibitively expensive. Drivers
will begin to
\emph{cruise for parking}~\cite{shoup:2005aa}, significantly contributing to surface street
congestion. Researchers have
attempted to measure the economic loss to both these drivers and the cities
themselves. For the former, drivers in different cities can spend anywhere
between 3.5 to 14 minutes searching for spots every time they
park~\cite{shoup:2006aa}. For the latter, cruising behaviors can lead to
substantial congestion in dense urban districts. For instance, there exists a
commonly cited folklore that 30\% of traffic in a city is directly due to
drivers looking for parking~\cite{shoup:2005aa}\footnote{It is not entirely clear where this number originates from as estimates vary wildly. A motivation of this paper is to rigorously quantify the actual effect of parking on congestion.}.

\begin{figure}[ht]
	\centering
	\includegraphics[width=0.4\textwidth]{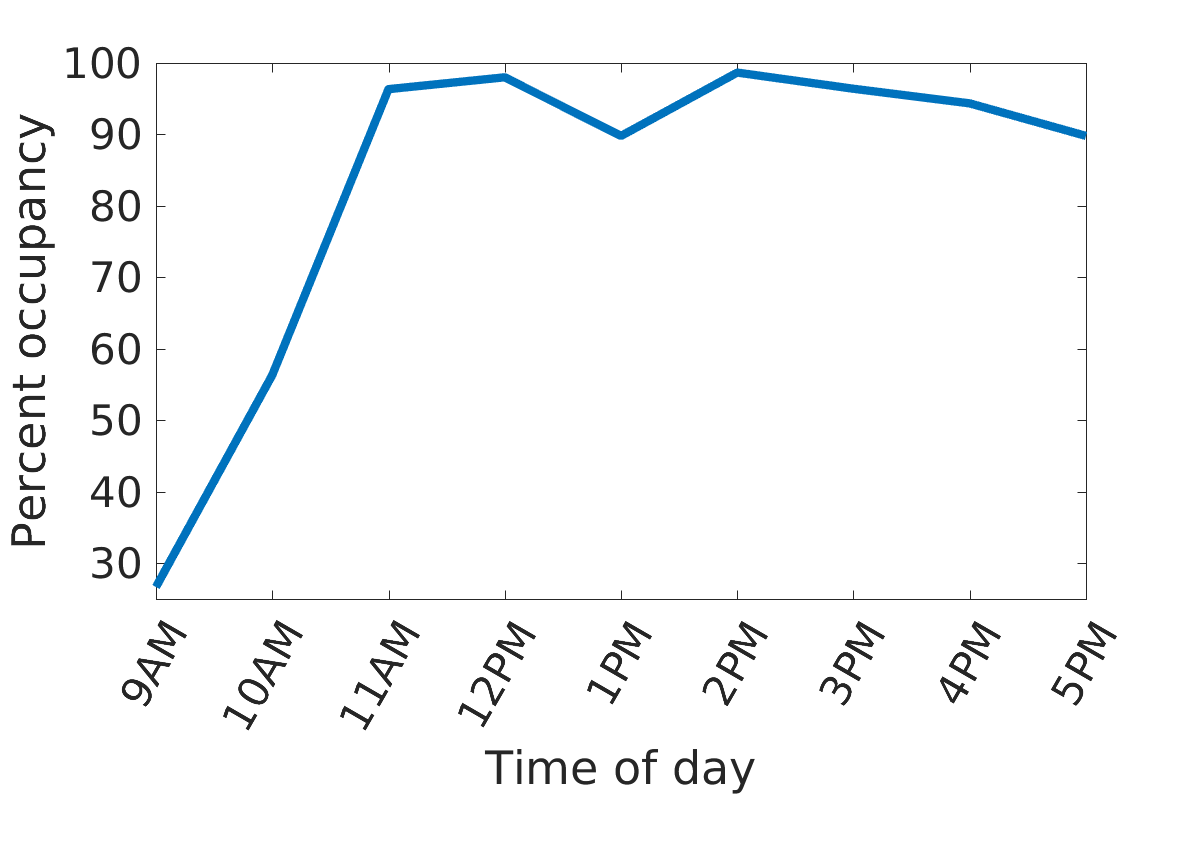}
	\caption{The observed parking occupancy (percentage of parking spaces in use), on a Saturday along the 3400 block of 18th St. in the Mission District of San Francisco. Fig. \ref{fig:meanvsblock} illustrates a key result of this paper: the congestion resulting from lack of parking along a block-face.}
	\label{fig:occup18}
\end{figure}

Municipalities and city planners typically aim to achieve some target \emph{occupancy}: the percentage of parking spaces in use at any given time~\cite{sfpark:2013aa}. Fig.~\ref{fig:occup18} shows the occupancy of the 3400 block of 18th St. in San Francisco, CA. Cities like San Francisco have launched projects like SFPark to target an average occupancy between around 85\% by slowly adjusting prices based on observed demand~\cite{sfpark:2013aa}.
\begin{figure}[!t]
	\centering
	\includegraphics[width=0.4\textwidth]{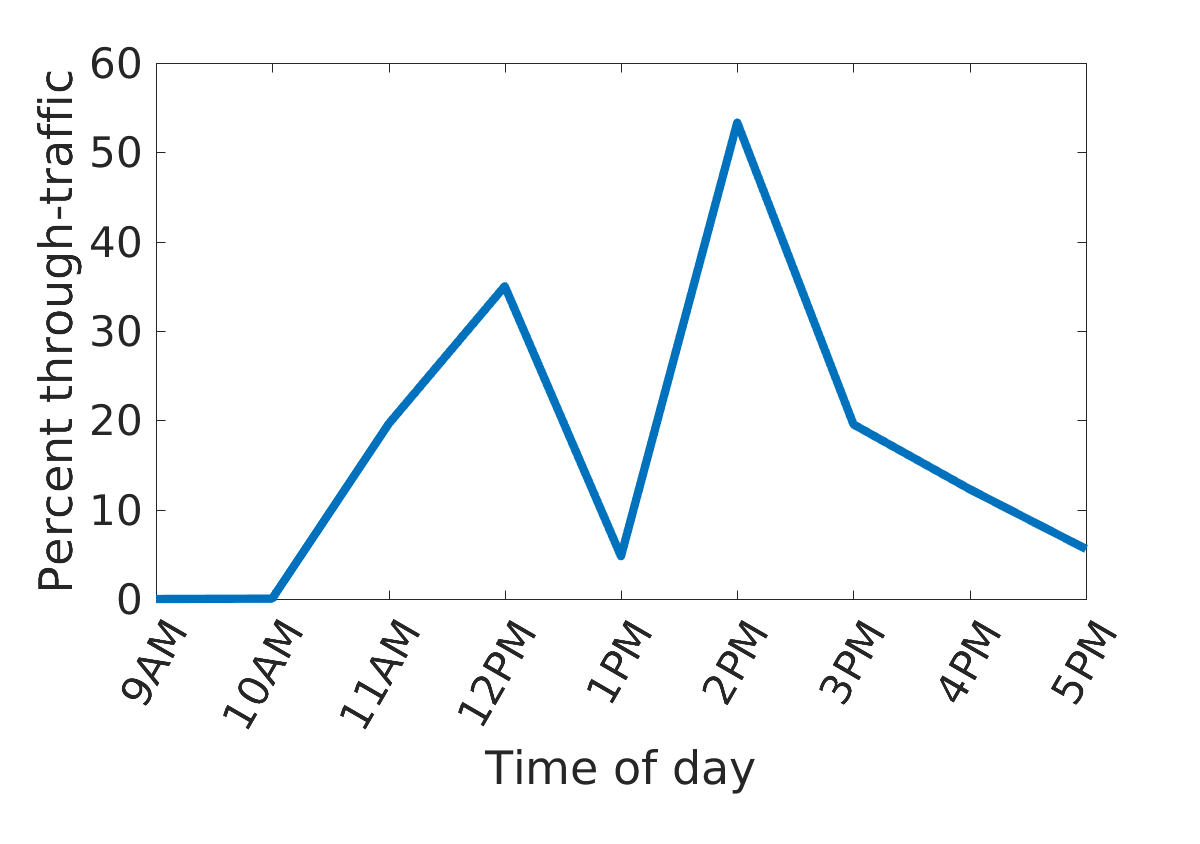}
	\caption{A visualization of one of our key results: estimated percentage of through traffic searching for parking. This estimate is obtained by determining the minimum arrival rate necessary to achieve an observed occupancy (Fig. \ref{fig:occup18}), and comparing this to the observed total through-traffic.}
	\label{fig:meanvsblock}
\end{figure}

Parking \emph{occupancy} (and availability) is an indirect measure (and means of control) of overall demand for vehicle access. Yet, if city planners must 
control congestion, occupancy alone is not a sufficient measure. Firstly, the
same occupancy levels of two streets in different parts of the city can lead to
different effects on through-traffic delays or respond differently to incremental price changes. Secondly, the street topology and interactions
between different blocks can lead to complex traffic dynamics, which a single
number like occupancy cannot capture. At the same time, cities cannot be overly
aggressive in controlling parking occupancy since they must maintain a high
availability of parking resources to serve downtown businesses and residents, as
well as delivery, courier, and emergency vehicle services. Therefore, a
reasonable question that a city planner would be interested in addressing is the following: \emph{Given a maximum tolerable level of congestion, what is the maximum occupancy at a block and what price achieves this occupancy?}


The question of parking's impact on congestion has remained difficult to address due to: 1) lack relevant data on pricing and demand and 2) lack of tractable and rigorous
models that link parking to congestion and capture spatial and temporal
variation. To address this question utilizing parking occupancy, traffic, and surface street topology data that is available today, our contributions are:
\begin{enumerate}
    \item \textit{Modeling}:~we describe and analyze a new kind of queue
	network where customers move between queues according to a network topology
	\emph{until} an available server is found, and leave the network after service
\item \emph{Control}:~we show that maximizing occupancy subject to
	constraints on the congestion created by drivers searching for parking is a convex program.
\item \textit{Application}:~we conduct a study based on real occupancy and pricing data for blocks in the San Francisco Mission District, showing that a) higher total occupancy does not necessarily lead to more traffic, and b) incentivizing drivers to park further away by reducing price can be equally as effective as disincentivizing drivers from parking at desirable locations.
\end{enumerate}

The paper is organized as follows. We provide motivation and review related work
in Section~\ref{sec:mot}. In Section~\ref{sec:model}, we present the network
queue model. We present results in Section~\ref{sec:results}. In particular, we
provide stability conditions under a uniformity assumption on the network
topology, we provide a
framework for determining the arrival rate in the non-uniform case, and we pose
an optimization problem to optimize parking availability subject to maximum congestion
constraints that we show to be convex. In Section~\ref{sec:application}, we
demonstrate the effectiveness of the solution to the optimization problem on a
network modeled after San Francisco's Mission District.
We conclude with discussion and commentary on future
work in Section~\ref{sec:conclusion}.

\section{Motivation}
\label{sec:mot}
As observed by Pierce and Shoup, circling for parking occurs when occupancy reaches
100\% \cite{pierce2013getting}, however, this takes an instantaneous point of view likely unavailable to city planners. Rather, if occupancy is
taken to be the expected proportion of parking spaces in use \emph{over a given time period}, then
high occupancy block-faces must be full at least some of the time, and therefore
responsible for some traffic---see, \emph{e.g.}, Fig.~\ref{fig:meanvsblock}. 

\subsection{Data Availability}
Municipalities (in particular, city planners and transportation departments) are gaining access
to data from recently installed \emph{smart parking meters} and, on occasion, individual parking space sensors (e.g., San Fransisco \cite{sfpark:2013aa}, Seattle \cite{ottosson:2013aa}, Los Angeles \cite{ghent2012express}, and Pittsburgh \cite{fabusuyi2013evaluation}). Yet, no city has completely implemented full-scale transportation sensor grids
that include active monitoring of parking on a space-by-space basis. Regardless
whether such a goal may be reached, however, many cities have a growing history
of parking transaction data collected by digital meters. These data can be used
to estimate parking occupancy; transactions provide an estimate of how long a driver intended to park and the number of drivers parked moment to moment. In our experiments, we make use of transaction, traffic, and infrastructural data publicly made available by the SFPark pilot study \cite{sfmta:2017aa}.



\subsection{Related Work}

Early work focused largely on parking supply and demand \cite{vickrey:1954aa}, and refinement of the economic view of parking continues through today \cite{inci:2015aa}. The costs of congestion caused by \emph{cruising for parking} \cite{shoup:2005aa, shoup:2006aa} have motivated research in modeling urban parking dynamics, and economizing of parking spaces has led to a desire to control demand levels via price.

Over the last few decades, a number of models (e.g., Vickrey's celebrated
``bathtub'' model) have been developed and introduced in the absence of data
only recently becoming available \cite{arnott2013bathtub, arnott:2006ab}. These
models typically take a time-varying flow and capacity view in the form of
systems of partial differential equations (see \cite{inci:2015aa} for an
overview of variations on these models).

Recent research has observed, however, that transaction data can be used to
estimate parking occupancy and, in consequence, used to estimate resource performance \cite{yang2017turning, dowling2017much}. The distinction that occupancy below 100\% results in congestion has recently been noted by~\cite{millard2014curb} in
their own analysis of the SFPark pilot study parallel to \cite{pierce2013getting}, however the authors of~\cite{millard2014curb} view block-face parking as a
Bernoulli random variable, between being full or not. We build on this work by 1) \emph{not} implicitly assuming curbside parking occupancy is independent between block-faces and 2) considering all possible states of parking spaces---from completely empty to completely full---along block-faces.

Occupancy and other data lend themselves to discrete and
probabilistic models that may potentially better reflect flow on surface streets
as compared to flow on highways or through spatially homogeneous regions, as in
\cite{arnott:2006aa} and \cite{arnott:2009aa}. Hence, classical methods of
queueing theory have recently been applied to parking areas: garage and curbside alike \cite{ratliff:2016aa, larson:2010aa, caliskan2007predicting, ceballos2004queue}. 

Our work primarily builds on existing parking literature by expressing curbside parking as a network of queues. Specifically, utilizing newly available parking data, we implement the basis for a spatially heterogeneous model city planners can use to effectively test parking policies and, furthermore, we determine that maximizing occupancy subject to congestion constraints using price controls is a convex optimization problem.

\section{Queueing Model}
\label{sec:model}

%
%
\subsection{Model Setup}


\begin{figure}[t]
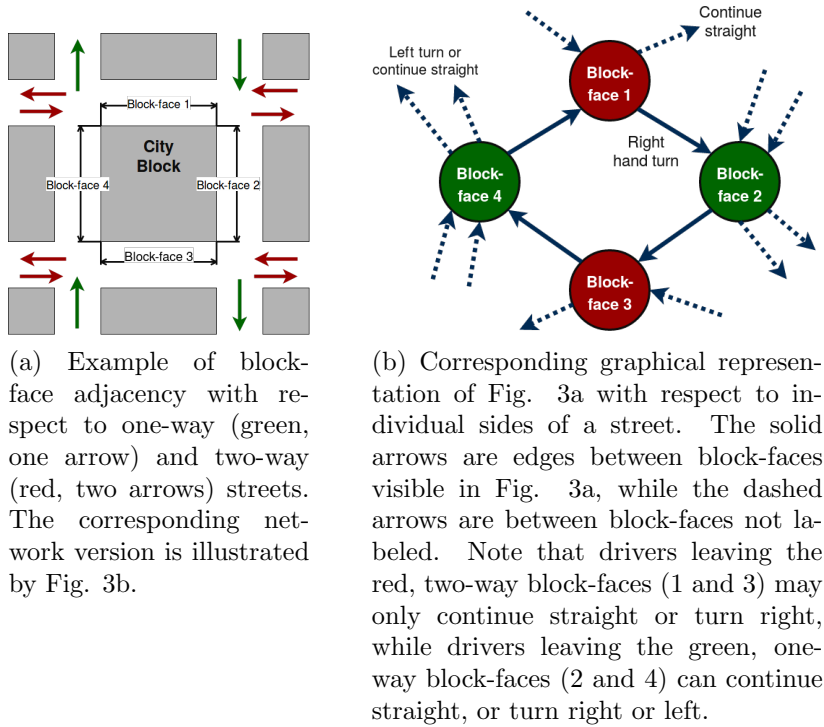

	\centering
	
	\subfloat[Example of blockface adjacency with respect to one-way (green, one arrow) and two-way (red, two arrows) streets. The corresponding network version is illustrated by Fig. \ref{fig:blocknetwork}.]{%
		\includegraphics[width=4cm]{figures/queueblock}%
		\label{fig:blockadj}%
	}\qquad
	\subfloat[Corresponding graphical representation of Fig. \ref{fig:blockadj} with respect to individual sides of a street. The solid arrows are edges between block-faces visible in Fig. \ref{fig:blockadj}, while the dashed arrows are between block-faces not labeled. Note that drivers leaving the red, two-way block-faces (1 and 3) may only continue straight or turn right, while drivers leaving the green, one-way block-faces (2 and 4) can continue straight, or turn right or left.]{%
		\includegraphics[width=6cm]{figures/queueblocknetwork}%
		\label{fig:blocknetwork}%
	}
	
	\caption{Block-face parking around a typical city block, and the corresponding graphical representation of the queue-network with respect to legal turns.}
\end{figure}


Although a natural model, queue networks have not been used extensively in parking related research~(see, e.g.~\cite{Raheja2010} and the references within for more details). Two major reasons for this are: 1) the size of the state space grows exponentially as the size of the network grows; 2) established queueing network results~(e.g., for communication networks) do not carry over directly. The rest of this section will describe the details of the queueing network model, its difference to conventional models, and how we overcome these difficulties.

\subsection{Queues Interacting Via Rejections}
We model each block-face as a multi-server queue, where the number of servers is
the number of available parking spots on that block-face. The block-faces are
connected as \emph{nodes on a graph}, where two nodes are adjacent if vehicles
can go from one block-face to the other in the road network. See
Fig.~\ref{fig:blockadj} for an example. To account for legal turning
maneuvers~(e.g., right turn only) and one way streets, we use \emph{directed
edges}. We use conventional notations $D =(V,E)$ to describe this digraph.
Without loss of generality (WLOG), we assume this graph is connected.

A queue, or a node $i \in V$ is characterized by an \emph{exogenous} arrival rate $\lambda_i$, a service rate
$\mu_i$, and the number of servers $k_i$. We assume that the \emph{exogenous} arrival process is Poisson~(independent between queues) and the services times are generally distributed like conventional $M/G/\cdot/\cdot$ queues~\cite{wolff1989stochastic}, however, unlike conventional queueing networks where customers are buffered at individual queues, we assume that customers (or drivers), are buffered or queued along the \emph{network edges}. This behavior reflects the key fact that vehicles which cannot find parking circulate in the network rather than wait at one location. Therefore, if the driver is served by a queue, it then leaves the network. However, if it finds the current queue to be full, it is \emph{rejected by that queue} and moves to neighboring queues to find new parking spots. The rate of these rejections is parking scarcity's contribution to through-traffic delays.

The key difference between our queue network and conventional networks--such as a Jackson network~\cite{jackson1957networks}--is that drivers proceed to other queues after they are rejected rather than served. Since the rejection of a queue with Poisson arrivals and exponential service times is not Poisson, characterizing the stationary distribution of this network of queues is very difficult because the distribution of total arrival rate itself to any queue is unknown.

Since the exact distribution of the queue is difficult to characterize, we instead turn to understanding the behavior of the \emph{mean performance metrics} of the network. This relaxation allows us to use theorems such as Little's Law~\cite{little1961proof} that do not depend on the exact distributions. Secondly, the controllable and measurable quantities are often average values like occupancy and parking service times.

\subsection{Stationary Distribution of a Single Queue}
Here we introduce how a \emph{single queue} can be analyzed, and later in the paper extend the analysis to a network of queues. To help avoid confusion between exogenous arrivals~(from outside of the network, denoted by $\lambda$) and endogenous arrivals~(rejection from neighboring queues, denoted by $x$), we use  $y$ as the total arrival rate to a queue. Suppose the service rate~(inverse length of parking time) of each server is $\frac{1}{\mu}$ and there are $k$ servers~($k$ parking spots) in total. Let $\pi_i$ be the stationary probability that $i$ servers are busy~($i$ cars are parked), for $i=0,\dots,k$. Let $\boldsymbol{\pi}=[\pi_0 \; \dots \pi_k]$. For this single queue, we can explicitly write down its stationary probability distribution via the transition rate matrix:
%

\begin{equation*}\label{eqn:transition}
   \bd Q= \bma -y & y & 0 & 0 & \cdots & 0 & 0 \\
   \mu & -(\mu+y) & y & 0 & \cdots & 0 & 0\\
   0 & 2\mu & -(2\mu+y) & y & \cdots & 0 & 0\\
   & & & & \vdots & & \\
   0 & 0 & 0 & 0 & \cdots & k\mu & -k\mu \ebma,
\end{equation*}
and $\boldsymbol{\pi}$ is the unique solution to
\begin{equation}\label{eqn:stationary}
\boldsymbol{\pi Q} = \boldsymbol{0}
\end{equation}
such that $\sum \pi_i =1$. Let $\rho = \frac{y}{\mu}$. By standard calculations~\cite{wolff1989stochastic},
\begin{equation}\label{eqn:statvalue}
\boldsymbol{\pi} = \pi_0 \cdot \left[ 1, \rho, \cdots, \frac{\rho^{k}}{k!} \right]
\end{equation}
where $\pi_0 = [\sum_{j = 0}^{k} \frac{\rho^{j}}{j!}]^{-1}$. Using Little's Law, the occupancy $u$, or the proportion of busy servers at any given time can be expressed as,
\begin{equation}\label{eqn:little}
u = \frac{y}{k\mu}\left( 1 - \pi_0 \frac{\rho^{k}}{k!} \right)
\end{equation}.

Note that $(1 - \pi_0 \frac{y^{k}}{k!})$ is the probability that \emph{at least} one space is available. Consider, if drivers are unable to wait for an available server at a particular block, in order to obtain occupancies approaching 100\%, cars would need to arrive at an infinite rate in order to immediately replace vehicles exiting service. Since it is often cited that congestion due to driver's searching for parking is a significant cost to the social welfare, this is a critical misconception.

A block-face queue is therefore rejecting incoming vehicles at a rate of  $y\cdot\pi_k$. The difficulty therein lies with estimating these total arrival rates, because no two adjacent block-faces are independent.

\section{Network of Queues}
\label{sec:results}
    In this section we study these networks of queues. We first consider the uniform case, then extend the results to the non-uniform case.


\subsection{Uniform Network}
    \label{sec:uniform}
Many urban centers have fairly uniform street topologies (e.g., the famed Manhattan streets), where the streets from a regular graph. In this section we make the assumption that the queueing network is entirely uniform: the topology is a $d$-regular graph, all block-faces have the same number of servers with the same service rate $\mu$, and they have the same exogenous arrival rate $\lambda$.


In this regular queue network, each queue will have equal stationary distributions in the steady state, therefore we only need to look at a single queue as representative of the state space of the entire network. Let $x$ be the average rate of rejection of a queue to one of its neighbors, and $dx$ be the total rejection to all of its neighbors. Let $y=\lambda+dx$ be the total arrival rate to a queue, where $\lambda$ is the exogenous arrivals and $dx$ are the rejections from its neighboring queues. We have the conservation equation,
\begin{equation}\label{eqn:consv}
    dx = y\pi_k,
\end{equation}
 where $\pi_k$ is the probability that all $k$ severs are busy.
 Combined with stationary distribution of \eqref{eqn:stationary} we have the following equations:
\begin{align}
	\left\{ \begin{array}{ll}
		\boldsymbol{\pi} \boldsymbol{Q} & = 0\\
		\sum \pi_i & = 1\\
		dx&=\pi_k (\la+dx)\end{array}\right.
\end{align}
We can write (\ref{eqn:consv}) as,
\begin{equation} \label{eqn:y}
y-\la= \frac{\frac{\rho^k}{k!}}{\sum_{i=0}^k \frac{\rho^i}{i!}} y
\end{equation}
where $\rho=\frac{y}{\mu}$. The equation in \eqref{eqn:y} is a polynomial in $y$. The next lemma states that there exists a unique solution to $y$ (and thus $x$) as long as the queues are stable:
\begin{lemma}
    \label{thm:lemma1}
 If $0 < \la < \mu k$, then (i) there is a unique and
    positive solution to $y$ in \eqref{eqn:y} and (ii) the solution is greater than $\la$. In addition, the rejection rate
    $x$ is also unique and positive.
    \label{lem:lem1}
\end{lemma}
The proof is given in Appendix~\ref{app:lem1}. This result states that as long as the total arrivals are less then the service rate times the number of spaces, we can explicitly find the rejection rates and the stationary probabilities by solving a polynomial equation. 
%

\subsection{Non-uniform Network}
    \label{sec:nonuniform}
Of course, the totally uniform assumption rarely holds up in practice. But given occupancy data we show that the \emph{total} exogenous and endogenous arrivals to a queue can still be solved for and used to estimate the traffic caused by drivers searching for parking. This time, for some \emph{total} incoming rejection rate $x$, letting $y = \lambda + x$, we can estimate the endogenous proportion of incoming arrivals as the sum of the outgoing fractional rejection rates of adjacent queues.

Assuming the queueing network reaches steady state,
from the perspective of a single queue in solving \ref{eqn:little} for $\pi_0$ gives
\begin{equation}
\pi_0 \frac{\rho^{k}}{k!}  + \frac{uk\mu}{y} = 1,
\end{equation}
where $u$ is the occupancy level and $\rho=\frac{y}{\mu}$.
Rearranging terms yeilds a polynomial in $y$,
\begin{equation}\label{eqn:y_poly_occup_mu}
0 = \sum_{i=0}^{k} \frac{1}{\mu^{i-1}}\left[\frac{i - uk}{i!}\right]y^{k}.
\end{equation}
Again, we can characterize the solutions to \eqref{eqn:y_poly_occup_mu}
\begin{lemma}\label{lem:lem2}
If $u \in [0,1)$ and $k$ is a positive integer, then \eqref{eqn:y_poly_occup_mu} has a unique real, positive root.
\end{lemma}
The proof is provided in Appendix~\ref{app:lem2}.

This root need not be bounded, hence the restriction of the values of $u$ to the interval $[0,1)$. In order to achieve a $100\%$ occupancy, implying the probability of being full is $1$, vehicles would need to arrive constantly ($y = \infty$), immediately taking the place of any vehicle that leaves upon completion of service. This is analogous to the requirement that for the $M/M/k/k$ queue to be stable, $\pi_0 > 0$.

\section{Optimizing parking availability}
    
\label{sec:optimization}

\label{sec:optproblem}

Price elasticity of demand provides a means of describing how consumer demand will change with incremental changes to price. Currently, Pierce and Shoup's analysis of the SFPark pilot project in \cite{pierce2013getting} is the state-of-the-art in estimating the price elasticity of demand for curbside parking; their exploratory analysis provided rough estimates of aggregated elasticities across time, location, and price change directions. For the purposes of this paper, and in order to make use of the results in \cite{pierce2013getting} we assume a \emph{linear} elasticity, however, any demonstrably reasonable (reflective of consumer behavior), concave function would not tax the validity of our results. Thus, a block-face has a linear elasticity
$\alpha$, and a function $\mathcal{U}: p \mapsto u$, taking a price $p$ to an occupancy level $u$, defined as

\begin{equation}\label{eqn:demand}
 \mathcal{U}(p) = 1 - \alpha p
\end{equation}

Recall \eqref{eqn:y_poly_occup_mu}; we can write
the right-hand side of this equation as a mapping $F: Y\times U\rightarrow
\mathbb{R}$ where $U=(0,1)$ such that
    \begin{equation}
        F(y,u)=\sum_{i=0}^{k} \frac{1}{\mu^{i-1}}\left[\frac{i - uk}{i!}\right]y^{k}
        \label{eq:new_y_poly}
    \end{equation}
    Note that this map is smooth in both its arguments $y$ and $u$.
   By applying the  Implicit
Function Theorem~\cite[Theorem C.40]{lee:2002aa}, a smooth mapping
$f:u\mapsto y$ exists and it is continuous and differentiable. Moreover, there
is an explicit expression for its derivative and the function $f$ maps an occupancy $u \in U$ to the unique real root
    $y$ of $F(y,u)=0$.

    Consider the following composition,
\begin{equation}\label{eqn:rejection}
 g(p) = f(\mathcal{U}(p)) \cdot \pi_k,
\end{equation}
which is equal the rate of rejection of vehicles from a block given a price $p$.
The composition \eqref{eqn:rejection} takes a price to a resulting level of
congestion along an edge in a queue network due to rejections.

The optimization problem given by
\begin{equation}\label{eqn:opt}
	\begin{aligned}
		& \underset{\boldsymbol{p}}{\text{maximize}}
		& & \mathcal{U}(\boldsymbol{p}) \\
		& \text{subject to}
		& & g_i(p_i) \leq \bar{x}_i, \; i = 1, \ldots, m.
    \end{aligned}\tag{P-1}
\end{equation}
maximizes parking resource utilization subject to a congestion constraints
$\bar{x}_i$ imposed on each block-face. Since \eqref{eqn:demand} is
concave, if $g_i$'s are convex, then \eqref{eqn:opt} is a convex optimization problem
easily solved by gradient descent.

\label{sec:convexity}
\begin{theorem}
    \label{lem:lem3}
	The optmization problem \eqref{eqn:opt} is convex.
\end{theorem}

The main technical challenge is to show that the constraints $g_i$'s are convex in the occupancy $u_i$. This is somewhat involved and the proof is provided in Appendix~\ref{app:lem3}.
Since \eqref{eqn:opt} is convex, there exist many ways to solve it. We use a projected gradient for the case study in the next section.



\section{APPLICATION}
\label{sec:application}
\begin{wrapfigure}{r}{0.5\textwidth}
	\vspace{-20pt}
	\begin{center}
		\includegraphics[width=0.3\textwidth]{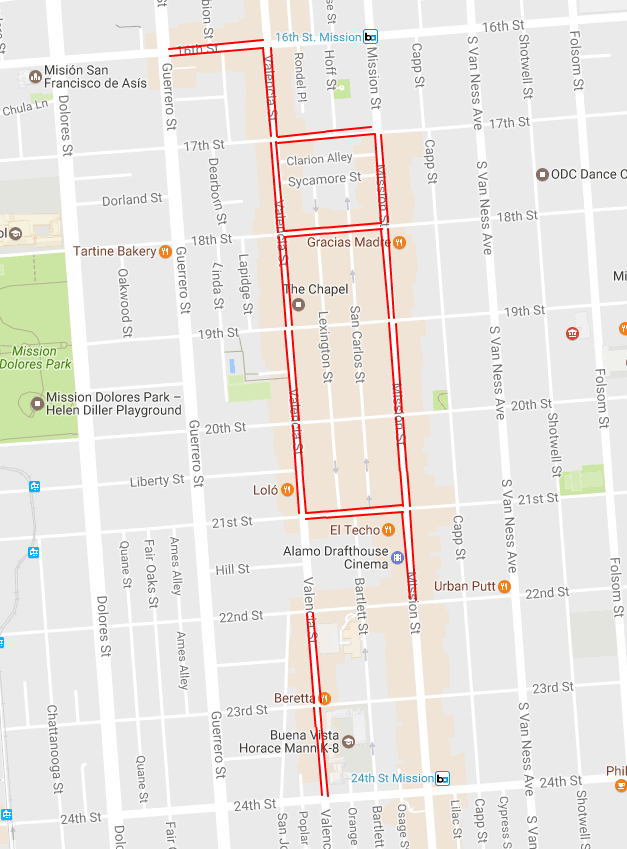}
		\caption{Block-faces, highlighted in red, with curbside parking data from the SFPark pilot program; Mission District, San Francisco.}
		\label{fig:mission}
	\end{center}
\end{wrapfigure}

    We consider the application of the above methods to curbside parking San Francisco's Mission District (Fig. \ref{fig:mission}). Using data collected by the SFPark pilot from May 8th, 2012 - August 29th, 2017 and elasticities estimated by \cite{pierce2013getting}, we identify block-faces responsible for the high congestion impacts to through-traffic and set constraints to bring this down to some hypothetically tolerable level. All data is calculated to an \emph{hourly} rate, e.g. the average percentage of parking spaces in use over the course of an hour.

According to \cite{pierce2013getting}, curbside parking in the Mission District
of San Francisco displayed an average price elasticity of $-0.21$. Price
elasticity varied greatly due to the time of day, week, and year, among a number
of other observable factors. For the purposes of demonstration in this paper, we assume a uniform price elasticity of $-0.21$ across block-faces in the Mission District, and therefore, resulting price changes should be taken with a grain of salt.

We examine two scenarios: 1) we wish to reduce overall \emph{congestion} due to parking by 80\% at two high occupancy block-faces and 2) achieve $>$80\% \emph{occupancy} at each block-face, rather than a neighborhoodwide average of 80\%, concentrated at a smaller proportion of the blocks in the district. We find that, in particular, \emph{spatial inefficiency, and not high occupancy, results in congestion.}

\subsubsection{Congestion Reduction}
\label{sec:expone}

\begin{wrapfigure}{r}{0.5\textwidth}
	\vspace{-20pt}
	\begin{center}
		\centering
		\includegraphics[width=0.5\textwidth]{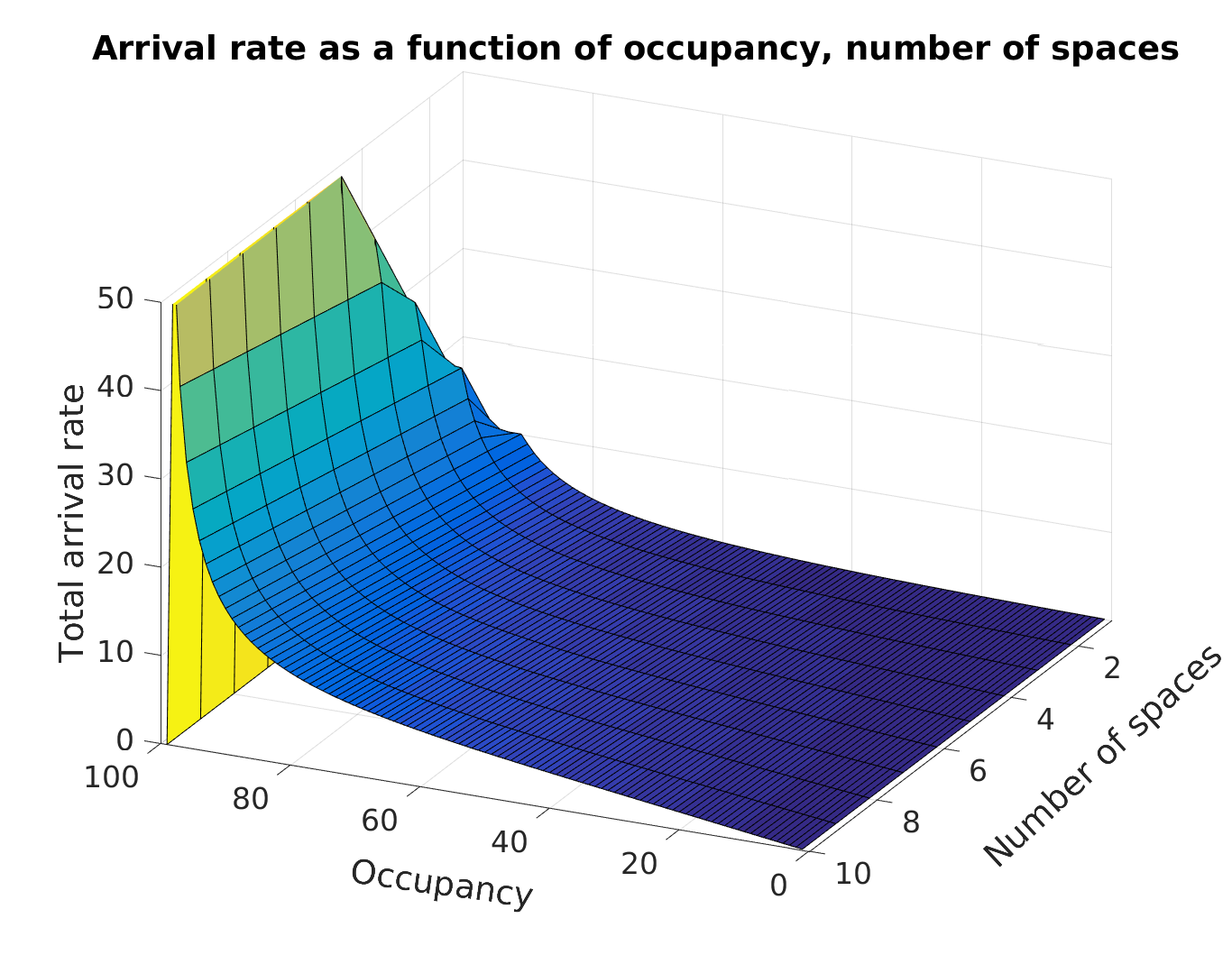}
		\caption{Neccesary total arrival rate $y$ to achieve an occupancy level for some fixed number of servers $k$ with a service time $\mu$ of 1. Not the sharp increase in total arrival rate around the 90\% occupancy mark and that increasing the number of servers only has a marginal bearing on this arrival rate.}
		\label{fig:y_as_occup_k}
	\end{center}
\end{wrapfigure}
\setlength{\belowcaptionskip}{-10pt}

The 3300 block of 17$^{th}$ street and the 3400 block of 18$^{th}$ street are responsible for the overwhelming majority of parking related congestion in Mission District at noon on the average Saturday, generating a total of nearly 60 vehicles unable to find parking per hour. As illsutrated by Fig. \ref{fig:meanvsblock}, a full third of 18$^{th}$ street's through traffic is made up of drivers unable to find parking.

At these traffic levels, 17$^{th}$ and 18$^{th}$ street have occupancies of 97\% and 98\% respectively. By increasing prices by \$0.28 on 17th$^{th}$ and \$0.27 on 18$^{th}$, we are able to reduce this congestion by 80\% to approximately 11 vehicles per hour, total, while still maintaining 91\% and 92\% occupancies, respectively. All other blocks see comparatively negligible changes.

The ``elbow'' of the highly non-linear curve describing the total arrival rate needed to achieve a particular occupancy level occurs around the 90\% mark, as illustrated in Fig. \ref{fig:y_as_occup_k}. By redistributing vehicles  intending to park at high occupancy blocks to historically low occupancy blocks through price control, less time is spent cruising for parking, leading us to our next experiment.

\subsubsection{Occupancy Redistribution}
\label{sec:exptwo}

\begin{wrapfigure}{r}{0.5\textwidth}
	\vspace{-20pt}
	\begin{center}
		\subfloat[Occupancy and resulting traffic in vehicles per hour generated.]{%
			\includegraphics[width=0.5\textwidth]{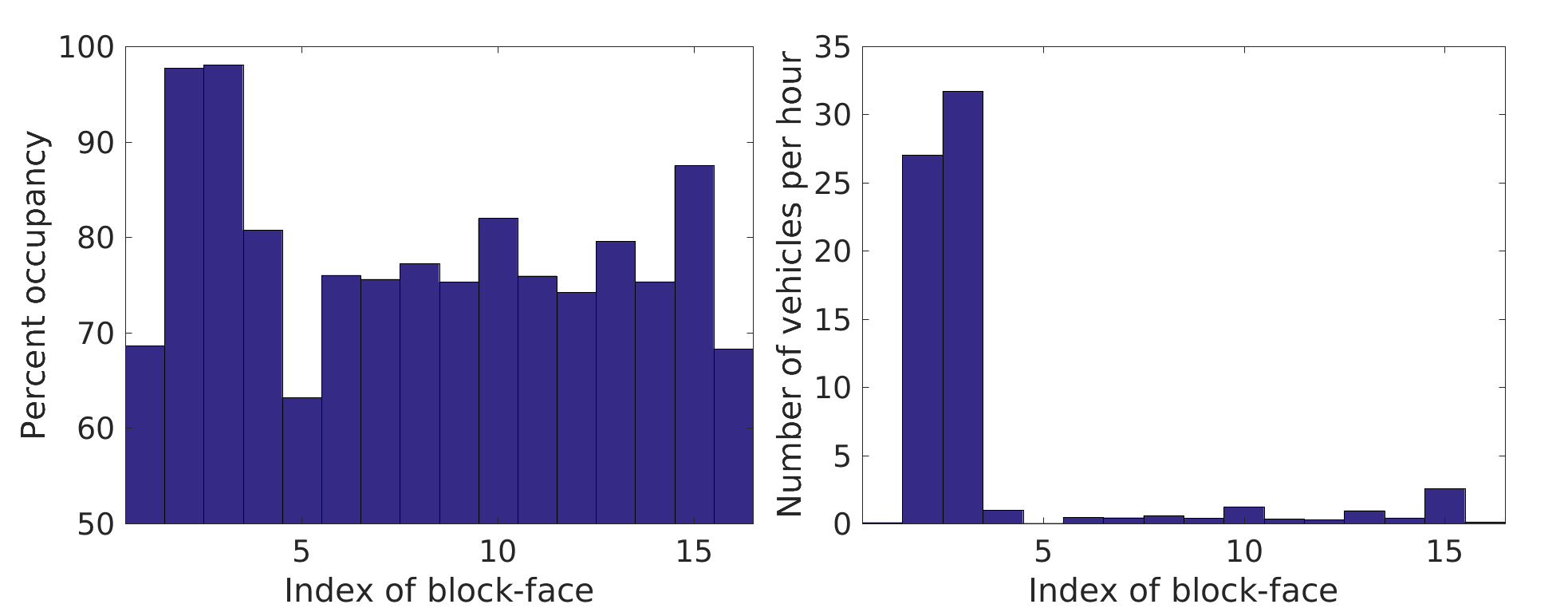}%
			\label{fig:uncontrolledparking}
		}
		
		\subfloat[Redistributing demand in Fig. \ref{fig:uncontrolledparking} to low-occupancy block-faces using the price changes indicated in Fig. \ref{fig:pricechange} results in less total traffic.]{%
			\includegraphics[width=0.5\textwidth]{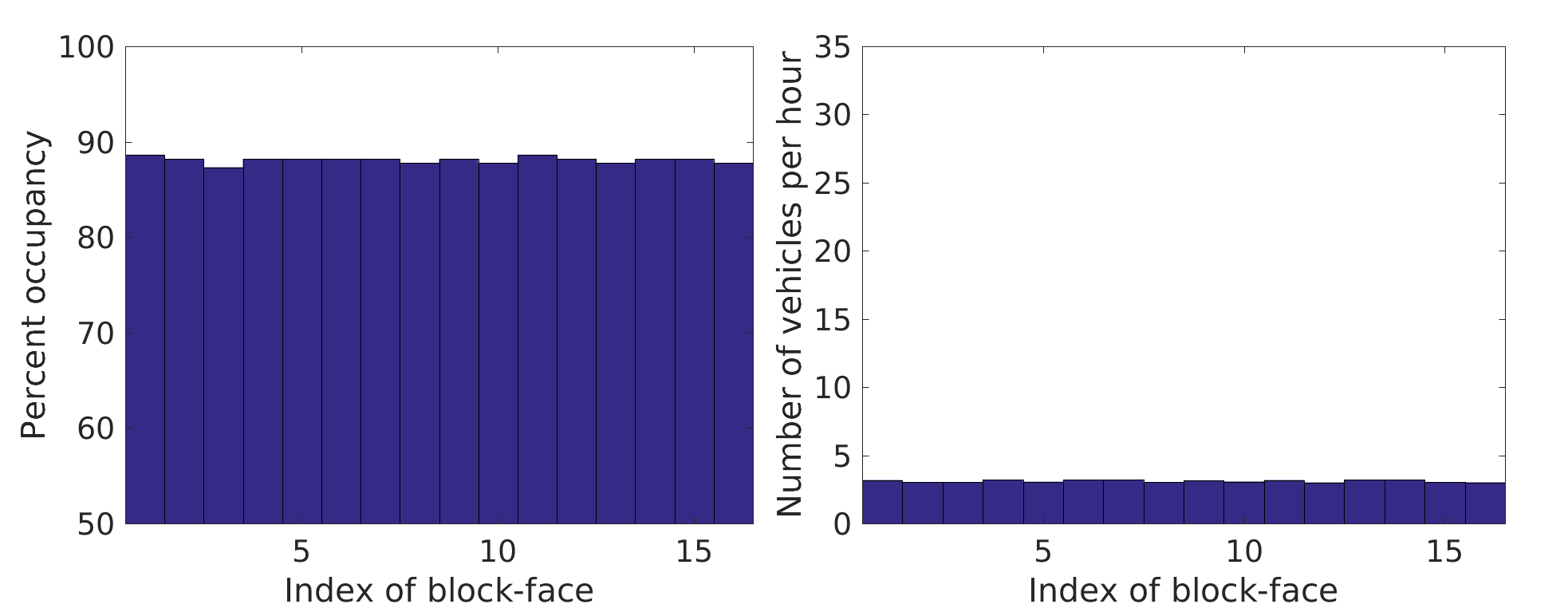}%
			\label{fig:controlledparking}
		}
		
		\subfloat[Price changes corresponding to the resulting occupancy redistribution in \ref{sec:exptwo}]{%
			\includegraphics[width=0.5\textwidth]{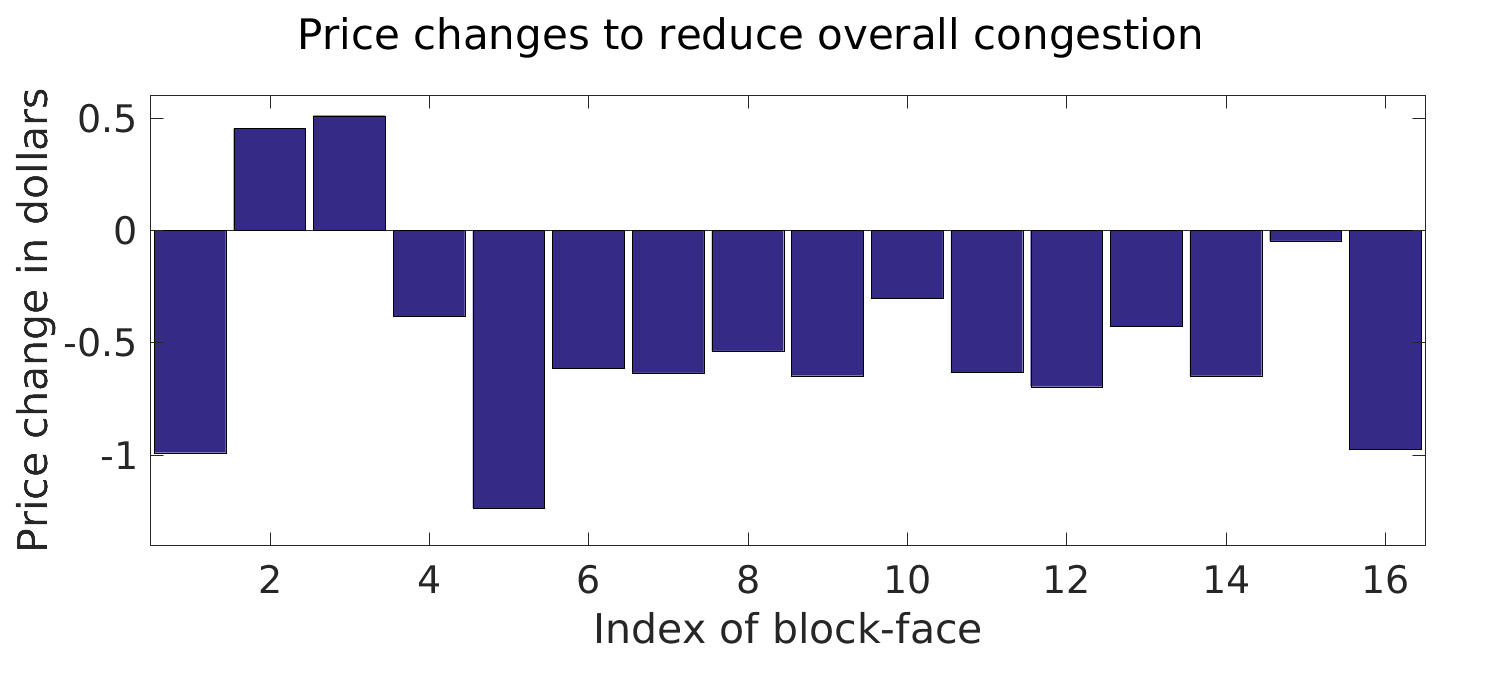}%
			\label{fig:pricechange}
		}
		
		\caption{Results of experiments in section \ref{sec:exptwo}}
		
	\end{center}
\end{wrapfigure}
\setlength{\belowcaptionskip}{-10pt}

On a typical Saturday at noon, the Mission District achieves an average occupancy of approximately 78\%, while generating over 60 vehicles per hour in additional traffic due to drivers searching for parking because there is a small number of high occupancy block-faces and a larger number of low occupancy block-faces. By bounding each block to producing no more than 1 vehicle every 20 minutes unable to find parking (for a total of 48 per hour for the district), each individual block-face individually exceeds 85\% occupancy \emph{at each block-face}. Indeed, after price control, the Mission District services a larger \emph{total} number of vehicles while still producing less additional traffic due parking scarcity.

Fig. \ref{fig:pricechange} indicates that, significantly discounting prices on low occupancy block-faces is an equally effective solution as raising prices at high occupancy block-faces, in order to achieve an effective distribution of parking resources that does not generate a costly amount of congestion searching for parking. Indeed, considering that a small number of block-faces may exhibit a high occupancy due to their desirable proximity to popular locations, incentivizing drivers to park somewhat further away may be more effective than pricing out other drivers by means of money or time to walk to a location.

\section{CONCLUSIONS AND FUTURE WORKS}
\label{sec:conclusion}
    \subsection{Conclusions}

With the growth of ride sharing services, electric vehicles, and increased demand for local delivery services, personal and commercial transportation is changing. In order for city planners to design effective future parking policies and make use of growing bodies of parking data, we developed a new kind of queueing network. We provided conditions for such networks to be stable, a ``single node'' view of a queue in such a network, and showed that maximizing the occupancy of such queues subject to constraints on the allowable congestion between queues searching for an available server is a convex optimization problem.

\subsection{Future Works}

A standing question in parking economics research is that of an appropriate maximum parking time \cite{inci:2015aa}. Some argue that a lower maximum parking time or lack of an initial buy-in price results in higher vehicle turn-over, and hence more congestion. Indeed, according to (\ref{eqn:little}), decreasing $\mu$ increases the total arrival rate necessary to achieve a fixed occupancy, but the probability of being full remains unchanged. Combined with the collection of ground-truth data and hypothesis testing, this question is closer to being answered.

Further, driver behavior is an important next-step to be considered. We have implicitly assumed that drivers, once inside the network searching for parking, will park regardless of price at a particular block-face. While this assumption alone is not unrealistic, how demand changes with respect to the total network sojourn time of the driver, distance from the initially desired location, and whether or not drivers have access to information regarding available parking locations are all certainly critical implications to consider.

\bibliographystyle{apalike}
		\bibliography{cdc2017}
%

\appendix
  \subsection{Proof of Lemma~\ref{lem:lem1}}

\label{app:lem1}
\begin{proof}
	Some algebra on \eqref{eqn:y} gives
	\begin{equation*}
\textstyle	k! (y-\la) \sum_{i=0}^k \frac{\rho^i}{i!}  = y\rho^{k}
	\end{equation*}
	The $\frac{y^{k+1}}{\mu^{k}}$ and $y\rho^{k}$ terms cancel, and we have a polynomial with degree $k$
	\begin{equation} \label{eqn:y_poly}
	\frac{\frac{k}{\mu^{k-1}}-\frac{\la}{\mu^{k}}}{k!} y^k + \frac{\frac{k-1}{\mu^{k-2}}-\frac{\la}{\mu^{k-1}}}{(k-1)!} y^{k-1} + \cdots + (1-\frac{\la}{\mu}) y - \la =0.
	\end{equation}
    Descartes' rule of signs~\cite{meserve:1982aa}, which roughly states that
    given a polynomial and ordering its terms from highest degree to lowest degree, the number of real positive roots is related to the number of sign changes. Let $n$ be the number of sign changes (from
    positive to negative), then the only possible number of positive roots to
    this polynomial are $n, n-2, n-4, \ldots$  In particular, if $n=1$, then the polynomial has one and only one positive root. Applying to the polynomial in \eqref{eqn:y_poly}, we notice the sign of the coefficients are determined by $k-\la$, $k-1-\la$, $k-2-\la$ and so on, until the constant term $-\la$. By assumption, $\la < \mu k$, so the first coefficient is positive. By assumption, $\la >0$, so the last coefficient (constant term) is negative. Then for any $\la \in (0,\mu k)$, it causes at most one change the signs of the other coefficients.
So $n=1$ for all possible $\la \in (0,k)$, and there is a unique positive solution to $y$.

To show that $y > \la$, let $f(y)$ be the polynomial in \eqref{eqn:y_poly}. We have $f(0)=-\la<0$, and $f(z) >0 $ for sufficiently large $z$ (positive coefficient on $y^k$ term). Since there is only one positive solution,  it suffices to show that at $f(\la)<0$. It turns out that $f(\la)$ has a telescoping sum, and
\begin{align*}
	f(\la) &= \textstyle \sum_{i=1}^k \frac{\la^i}{(i-1)!} - \sum_{i=1}^k \frac{\la^{i+1}}{i!} - \la \\
	&= \textstyle  \la - \frac{\la^{k+1}}{k!}-\la \\
	& < 0.
\end{align*}\end{proof}

\subsection{Proof of Lemma~\ref{lem:lem2}}
\label{app:lem2}



\begin{proof}
Let us first examine the coefficients of  $y^{k}$. WLOG, assume $\mu = 1$. We have the following sequence:
\begin{equation}\label{eqn:coeff}
s =\textstyle \{-uk, 1 - uk, \frac{ 2 - uk}{2!}, \ldots , \frac{k - uk}{k!} \}
\end{equation}
We will show that if $u \in [0,1)$, $k \in \mathbb{Z_{+}}$, the sequence
    (\ref{eqn:coeff}) undergoes exactly 1 sign change, and again apply Descartes'
    rule of signs. Observe that $s_0 < 0$ for any allowable values of $u$ and
    $k$. Further, observe that $s_k=(1-u)\left( (k-1)! \right)^{-1}$.
By induction, $s_k$ will always be positive for any value of $k$. If $k = 1$,
then $s_1 = (1 - u)(1)^{-1}$, and since $u \in [0,1)$, $s_1 > 0$. Assume this is
    true for $k$, then for $k + 1$, $s_k=(1-u)\left( k! \right)^{-1}$,
%
so that we have that $s_{k+1} > 0$. It now suffices to show that $\{s\}$ can
only undergo one sign change as we increment $i$. For some $k$, the $i$--th
element of $\{s\}$ is $s_i=(i-uk)(i!)^{-1}$.
%
Fix $k$. While the denominator of the sequence is itself increasing with $i$
(meaning $\{s\}$ need not be monotonic), it is strictly positive. We need only
look at the sign of the numerator. In particular, $uk$ is fixed between
$[0,1)\cdot k = [0,k)$, and $i$ is the set of indices between $[0,k]$. The
sequence (\ref{eqn:coeff}) will be negative until $i > \floor{uk}$, and since
$\floor{uk} < uk$, we are ensured there is only one sign change.

Since the coefficients of \eqref{eqn:y_poly_occup_mu} undergo one sign change, we again invoke Descartes' rule, and observe that we have one real positive root.
\end{proof}

\subsection{Proof of Theorem~\ref{lem:lem3}}
\label{app:lem3}
\begin{proof}
	Let $x = ku$. Then we can think of \eqref{eqn:y_poly_occup_mu} as
	\begin{equation}\label{eqn:xypoly}
	F(y,x) = \textstyle (\frac{x}{k!} - \frac{1}{(k-1)!})y^{k} + \cdots +
    (\frac{x}{2!} - 1)y^{2} + (x - 1)y + x
	\end{equation}
	Implicit differentiation of \eqref{eqn:xypoly}, written as $D_xF+ D_yF\cdot
    y'$ where $y'=dy/dx$,
    gives
	\begin{align}
	0&=\textstyle(\frac{y^{k}}{k!} + \cdots + y + 1) + ((\frac{1}{(k-1)!} -
    \frac{x}{k!})ky^{k-1} +\notag\\
    &\ \ \textstyle \cdots + (1 - x))y'
	\end{align}
Noting that  $(D_xF)(y)=\frac{y^{k}}{k!} + \cdots + y + 1$ and $(D_yF)(x,y)=(\frac{1}{(k-1)!} -
    \frac{x}{k!})ky^{k-1}+\cdots+(1-x)$ so that
    \begin{equation} 
        y'=-D_xF\cdot (D_yF)^{-1}
    \label{<++>}
\end{equation}
\begin{prop} \label{prop:y}
	Let $(x,y)$ be a positive solution to $F(x,y)=0$, then $y'$ evaluated at that solution is positive.
\end{prop}
We first show the theorem assuming the proposition is true.
We can similarly compute the second order implicit derivative $d^2y/dx^2$;
indeed,
\begin{equation}
    y''=\frac{D_xF\cdot( D_y^2F\cdot y'+
    D_{x,y}F)-D_yF\cdot D_{y,x}F\cdot y'}{(D_yF)^2}
    \label{eq:ypp}
\end{equation}
Hence, if $D_xF\cdot( D_y^2F\cdot y'+
    D_{x,y}F)-D_yF\cdot D_{y,x}F\cdot y'>0$ then $y''>0$. We have
    \begin{align}
       & D_xF\cdot( D_y^2F\cdot (-D_xF\cdot (D_yF)^{-1})+\\
    &D_{x,y}F)-D_yF\cdot D_{y,x}F\cdot (-D_xF\cdot (D_yF)^{-1})\notag\\
    &=D_xF\cdot (D_y^2F\cdot (-D_xF\cdot (D_yF)^{-1})+2 D_{y,x}F)\\
    &=D_xF\cdot {h}(x,y)
        \label{eq:deriv}
    \end{align}
    where ${h}(x,y)=D_y^2F\cdot y'+2D_{y,x}F$. Since $D_xF >0$, we focus on $h(x,y)$:
    Now,
    \begin{equation}
        (D_{y,x} F)(y)=((k-1)!)^{-1}y^{k - 1} + \cdots + 1
        \label{eq:dyx}
    \end{equation}
    and
    \begin{equation}
        -D_y^2 F =\textstyle(\frac{x}{k!} - \frac{1}{(k - 1)!}) k(k-1)y^{k -2} +
        \cdots +2(\frac{x}{2} - 1)
        \label{eq:yy}
    \end{equation}
Collecting all the $x$ terms in $D_y^2 F$ we can define
\begin{equation}
    \tilde{h}(x,y)=\textstyle\frac{x}{(k-2)!}y^{k -2}+\cdots+x.
    \label{eq:xterms}
\end{equation}
Since $F(y,x)=0$, we have
\begin{align}
    \textstyle \frac{x}{k!}y^k+\frac{x}{(k-1)!}y^{k-1}+\cdots+x
    =\frac{1}{(k-1)!}y^k+\cdots +y
    \label{eq:yterms}
\end{align}
so that
\begin{align*}
    &\textstyle \tilde{h}(x,y)+\frac{x}{k!}y^k+\frac{x}{(k-1)!}y^{k-1}- \frac{x}{k!}y^k-\frac{x}{(k-1)!}y^{k-1}\\
    &\textstyle
    =\frac{1}{(k-1)!}y^k+\cdots+ y\textstyle-\frac{x}{k!}y^k-\frac{x}{(k-1)!}y^{k-1}
    \label{eq:xterms2}
\end{align*}
Then,
\begin{align*}
 D_y^2 F &=\textstyle\frac{x}{k!}y^k+\frac{x}{(k-1)!}y^{k-1}+\frac{k}{(k-2)!}y^{k-2}+\cdots+2\\
 &-\textstyle\frac{1}{(k-1)!}y^k-\cdots - y.
\end{align*}
so that
\begin{align*}
h(x,y)&=\textstyle\frac{2}{(k-1)!}y^{k - 1} + \cdots + 2-\left(\textstyle\frac{1}{(k-1)!}y^k+\cdots
   + y\textstyle-\frac{x}{k!}y^k \right.\\
 &\left. \textstyle -\frac{x}{(k-1)!}y^{k-1}  -\frac{k}{(k-2)!}y^{k-2}-\cdots-2\right)y' \\
 & = \textstyle y'\left(\frac{x}{k!}-\frac{1}{(k-1)!}\right) y^k \\
 & \textstyle + y'\left(\frac{2}{(k-1)!y'}+\frac{x}{(k-1)!}-\frac{1}{(k-2)!}\right) y^{k-1} \\
 & \textstyle + y'\left(\frac{2}{(k-2)!y'}+\frac{k}{(k-2)!}-\frac{1}{(k-3)!}\right) y^{k-2} \\
  & \textstyle + y'\left(\frac{2}{(k-3)!y'}+\frac{k-1}{(k-3)!}-\frac{1}{(k-4)!}\right) y^{k-3} \\
 & \vdots \\
 & + \textstyle y'\left(\frac{2}{y'}+2\right).
 \end{align*}
 Through straightforward, but somewhat cumbersome algebra, we can show that if $(x,y)$ is a pair such that $F(x,y)=0$, then
 \[ \frac{2}{y'}+1 \geq x. \]
 Following the above inequalities and using $\frac{2}{y'}+2 \geq x$, at the solution $(x,y)$ where $F(x,y)=0$
 \begin{align*}
	 h(x,y) & \geq \textstyle y'\left(\frac{x}{k!}-\frac{1}{(k-1)!}\right) y^k \\
	 & \textstyle + y'\left(\frac{x}{(k-1)!}-\frac{1}{(k-2)!}\right) y^{k-1} \\
	 & \textstyle + y'\left(\frac{x}{(k-2)!}-\frac{1}{(k-3)!}\right) y^{k-2} \\
	 & \vdots \\
	 & + \textstyle y'\left(x\right) \\
	 &= y' F(x,y) \\
	 &=0,
	 \end{align*}
	 and $y'' \geq 0$ follows from $h(x,y)\geq 0$.

Now we prove Prop. \ref{prop:y}. This lemma follows from the Gauss-Lucas Theorem~\cite{Meserve1982}, which states that if $p(z)$ is a polynomial with real coefficients with complex roots $r_1,\dots,r_n$, then the complex roots of $p'(z)$ is contained in the convex hull of $r_1,\dots,r_n$. For a fix $x$, applying this theorem to $D_yF$ yields the fact that real parts of all roots of $D_yF$ is less than the root of $F(x,y)$. Since $D_yF \rightarrow -\infty$ as $y\rightarrow \infty$, at the root of $F(x,y)$, $D_yF \leq 0$. By \eqref{<++>} and the fact $D_xF >0$, $y'>0$. \end{proof}

\end{document}